\documentclass[11pt,onecolumn]{article}
\usepackage{amscd}
\usepackage{times}
\usepackage{amsmath}
\usepackage{amssymb}
\usepackage{xspace}
\usepackage{theorem}
\usepackage{graphicx}
\usepackage{ifpdf}
\usepackage{url,hyperref}
\usepackage{latexsym}
\usepackage{euscript}
\usepackage{xspace}
\usepackage[all]{xy}
\usepackage{color}
\usepackage{makeidx}
\usepackage{tikz}
\long\def\remove#1{}
\newtheorem{theorem}{Theorem}[section] 

\newtheorem{obs}[theorem]{Observation}
\newtheorem{corollary}[theorem]{Corollary}

\newtheorem{definition}[theorem]{Definition}
\newtheorem{proposition}[theorem]{Proposition}

\newcommand {\mm}[1] {\ifmmode{#1}\else{\mbox{\(#1\)}}\fi}

\newcommand{\supp}{\mathrm supp}
\newcommand{\coker} {\mathrm coker}

\newcommand{\cancel}[1]

\begin{document}

\title {Barcodes for closed one form - an alternative to Novikov theory} 

\author{
Dan Burghelea  \thanks{
Department of Mathematics,
The Ohio State University, Columbus, OH 43210,USA.
Email: {\tt burghele@math.ohio-state.edu}}
}
\date{}
\maketitle

\begin{abstract}

We extend the configurations $\delta_r$ and $\gamma_r$, discussed in \cite {B} and \cite{BH}, equivalently the closed, open and  closed-open   
bar codes from real- or angle-valued maps,  to topological closed one forms on compact ANRs. As a consequence one provides an extension of the classical Novikov complex associated to a closed smooth one form and a vector field the form is Lyapunov for,  to a considerably larger class of situations.  
We establish strong stability properties and Poincar\'e duality properties when the underlying space  is a closed manifold.  Applications to Geometry, Dynamics and Data Analysis are the targets of our research. A different approach towards such bar codes was proposed in Usher-Zhang's work
cf. \cite{UZ}.

\end{abstract} 

This paper is essentially my lecture at the workshop "Topological data analysis meets symplectic topology", Tel Aviv,  April 29- May 3, 2018.
 A beamer version can be downloaded from 
\newline https://people.math.osu.edu/burghelea.1/notes/index.html.  A detailed version of this work will be posted shortly.

\section {Classical Novikov theory}

Classical Novikov theory considers

- a smooth manifold $M^n$ (for the purpose of this paper a closed manifold), 

- a pair $(X,\omega), X $ a smooth Morse-Smale vector field on $M,$ $\omega\in \Omega^1(M)$ a smooth closed one form Lyapunov for $X,$

- a field $\kappa$ 

\noindent and relates 

-   the rest points  of $X$  and

-   the instantons (isolated trajectories) between  rest points of $X$  

\noindent to the topology of $(X; [\omega])$ \footnote{ $[\omega]$denotes the cohomology class represented by $\omega$} and the zeros of $\omega.$

The Morse-Smale property of $X$ implies that 

-   each rest point of $X$ is hyperbolic and has a {\it Morse index} $r,$ an integer between $0$ and $n,$ and 

-   all {\it isolated trajectories}, i.e. instantons between rest points, exist only  from rest points of index $r$ to rest points of index $(r-1).$  

The Lyapunov property of $\omega$ identifies the rest points of $X$ to the zeros of $\omega.$  
\vskip .1in

\noindent The theory constructs:

-  a  field $\kappa(\omega),$ extension of $\kappa,$ called  the {\it Novikov field} associated to $\omega,$

-  a  chain complex of finite dimensional vector spaces over $\kappa(\omega)$ whose  

\hskip .1in  a)  {\it components  $C_r$} 
are the vector spaces generated by the rest points of index $r$, and

\hskip .1in b) {\it  boundary maps  $\partial _r:C_r\to C_{r-1}$} 
are given by matrices whose entries are derived from the instantons  between  rest points 
appropriately counted, called "instanton-counting matrices".  The result of the "appropriate counting",  even when the cardinality of the instantons is infinite, is an element in $\kappa(\omega).$ The rank $\rho_r$ of  $\partial _r$  indicates that  there are at least $\rho_r$ different rest points of index $r$  which are the origin of instantons  and at least $\rho_r$ different rest points  of index $r-1$ which are the end of instantons. 

The main result of the Novikov theory is that the dimensions of the homologies of this complex are the Novikov-Betti numbers $\beta^N_r(X;[\omega])$  
and that the complex, up to a non canonical isomorphism, depends only on the closed one form $\omega.$  Elementary linear algebra arguments imply that the rank  $\rho_r$ plus the rank $\rho_{r-1}$ is equal to  the number $c_r $ of rest points of $X$ of index $r$ minus  the Novikov-Betti number $\beta^N_r(X;\omega),$ 
hence the complex, up to an isomorphism, is determined by either two of these three  types of numbers: $\beta^N_r, \  \rho _r, \  c_r.$   
 
Recall that any chain complex of vector spaces over a field $\kappa$ is non canonically isomorphic to a complex of the following type :
$$\xymatrix{ \cdots \ar[r] &(C_r= \mathcal H_r \oplus C_r ^+ \oplus C_{r}^-) \ar[r]^-{\partial _r}& (C_{r-1}= \mathcal H_{r-1} \oplus C_{r-1} ^+ \oplus C_{r-1}) \ar[r]\ar[r]
^-{\partial _{r-1}}&\cdots}$$
with $C_{r}^-  = C_{r-1}^+,$ and $$\partial _r= \begin{pmatrix} 0&0&0\\0&0&Id\\0&0&0\end{pmatrix},$$ called Hodge type.
 
 The chain complex $(C_\ast,\partial _\ast)$
 
-  {\it provides the exact number of rest points of Morse index $r$},

- {\it give a good information about the instantons between rest points},

- {\it establishes the existence of closed trajectories}
(when combined with standard Betti numbers),

\noindent results of interest even outside mathematics if the problem under consideration can be modeled by the input of Novikov theory.

The drawbacks of the theory for applications outside mathematics  are :

\begin{enumerate}
\item  reduced generality (manifold structure for the underlying space, Lyapunov closed smooth one form which is Morse used to represent  the "action" controlling  the  dynamics defined  by the vector field, 
\item  the infinite cardinality of the mathematical objects  involved in the definition of  the field $\kappa(\omega)$ and of the matrix $\partial _r$ which makes these objects not computer friendly \footnote{ in this paper the adjective {\it computer friendly} refers to concepts consistent with or invariants computed by 
 a computer implementable algorithm}.
\end{enumerate}

In view of potential interest outside mathematics 
we like:     

\vskip .1in
\begin{enumerate}
\item to derive the collection of numbers $ \beta_r, \ \rho_r,\  c_r$ as {\it computer friendly invariants} and without any reference to the  field $\kappa(\omega),$ 

\item to extend the result to  a  larger class of situations, i.e. a compact ANR denoted from now on by $X$  
instead of manifold $M^n$ , topological closed one form $\omega$ on $X$ (to be defined in section  2) instead of a smooth closed one form $\omega$ on $M^n.$
\end{enumerate}

This was already done in \cite{B} and implicitly in \cite{BH} in the case the topological closed one form  $\omega$ is exact and in the case $\omega$ represents an integral cohomology class.  These two situations represent the case of topological closed one form of degree of irrationality $0$ and $1.$ In the work summarized in this paper  we go much further and we treat the case of an arbitrary  {\it topological closed one form} from now  on abbreviated to TC1-form of arbitrary degree of irrationality.

In the theory referred to as AMN (Alternative to Morse-Novikov theory) the underlying space is always a compact ANR (absolute neighborhood retract), cf. \cite {B} chapter 1 for definition.
All spaces homeomorphic to finite simplicial complexes and in particular the compact  manifolds, or  more general the compact stratified spaces,  are  compact ANRs. 

The concept of tame TC1-form, a replacement (generalization) of a Morse closed one form,  is recalled in section 2,  the set $\mathcal O(\omega),$ {\it orbits of critical values} associated to a TC1-form  $\omega,$ are introduced in section 3 and the main results  are formulated in section 4 below.  

 Note that in \cite {B} to a tame integral TC1-form, equivalently an angle-valued map, $\omega$ one associates a principal $\mathbb Z-$covering $\tilde X\to X$ and  lifts $f:\tilde X\to \mathbb R$  of $\omega.$  These are continuous maps which are proper, with the homology of the levels $f^{-1}(t)$ finite dimensional  vector spaces and with the set of critical values  discrete subsets of $\mathbb R.$ All these properties were used in \cite{B}.

In this paper to  an arbitrary $\omega$ one associates the principal $\mathbb Z^k-$covering $\tilde X\to X$ and  lifts  $f:\tilde X\to \mathbb R$  but in the case the degree of irrationality (the number $k$) is $>1,$ 

-  the map $ f$ is never proper, 

-  the  levels $f^{-1}(t)$ might not have the homology as finite dimensional vector spaces,   

- the set of critical values if  nonempty is never  discrete  but  always dense.  

\noindent The approach of level persistence  considered in \cite{B} chapter 4, via graph representations, is apparently not applicable. 

However, the treatment described in \cite{B} sections 5 and 6, can be refined and leads to barcodes which are recorded as real numbers \footnote {this because a closed one form is view as a real-valued function up to a translation, cf definition 2 in section 2, so the barcodes of $f$ are intervals (i.e. barcodes of real-valued functions) up to a translation},  namely as configurations of points in $\mathbb R$  and $\mathbb R_+$ rather than $\mathbb R^2$ or $\mathbb R^2/\mathbb Z.$  

\section{Topological closed one forms (TC1-forms) and tameness}\

A topological closed one form (TC1-form) extends the concept of  smooth closed one form on a smooth manifold  to an arbitrary topological space. One way to obtain this is  to view it  as an equivalence class 
of multivalued maps (first definition) an other way is as an equivalence class of equivariant maps on an associated principal $\mathbb Z^k-$covering (second definition). 
\vskip .2in

{\bf First definition}
\begin{definition}\

\begin{enumerate}
\item A {\bf multi-valued map} is  
a systems $\{U_\alpha, f_\alpha: U_\alpha\to \mathbb R, \alpha\in A\}$ s.t. 
\hskip .1in \begin{enumerate}
\item $U_\alpha$ are open sets with $X=\cup U_\alpha,$   
\item $f_\alpha$ are continuous maps s.t $f_\alpha- f_\beta: U_\alpha\cap U_\beta\to \mathbb R$ is locally constant 
\end{enumerate}
\vskip .1in
\item 
Two  multi-valued maps are {\bf equivalent} if together remain a multi-valued map.
\item  A  {\bf TC1-form}   
is an equivalence class of multi-valued maps. 
\vskip .1in
\end{enumerate}
\end{definition}

{\it Examples:} 
\begin {enumerate}
\item A smooth  closed one form $\omega\in \Omega^1(M)$ with $d\omega=0$ defines  a TC1-form  $\omega$. 

Indeed in view of Poincar\'e Lemma for any $x\in M$ choose an open neighborhood  $U_x\ni x$  and $f_x:U_x\to \mathbb R$ a smooth map s.t. $\omega_x|_{U_x}= d f_x$. The system $\{U_x, f_x:U_x\to \mathbb R\}$ provides a representative of the TC1-form $\omega$.
\item A simplicial one cocycle $\omega$ on the simplicial complex $X$ defines  a TC1-form.

Indeed, if $X$ is a simplicial complex, $\mathcal X_0$ the collection of vertices and $S\subset \mathcal X_0\times \mathcal X_0$ consists of pairs $(x,y)\in \mathcal X_0$ s.t. $x,y$ are the boundaries of a $1-$simplex  of $X$ then a simplicial one cocycle is a map $\delta: S\to \mathbb R$ with the properties $\delta(x,y)= -\delta(y,x)$ and for any three vertices $x,y,z$ with $(x,y), (y,z), (x,z)\in S$ one has $\delta(x,y)+ \delta(y,z) + \delta(z,x)=0.$
The collection of open sets  $U_x,$ $U_x$ the open star of the vertex $x,$ and the maps $f_x: U_x\to \mathbb R,$ the linear extension on each open simplex of $U_x$ of the map defined by $f_x(y)= \delta(x,y)$ and $f_x(x)=0,$ provide a representative of the TC1-form defined by the one cocycle $\delta.$ 
\end{enumerate}
\vskip .1in

Clearly a TC1-form $\omega$ determines a cohomology class $\xi(\omega)\in H^1(X;\mathbb R).$
\footnote{
It suffices  to show that a representative $\{ U_\alpha, f_\alpha: U_\alpha\to \mathbb R\}$ of $\omega$ defines, for any closed continuous path $\gamma: [a,b]\to X,$  the number $\int_{\gamma} \omega\to \mathbb R$ independent on the homotopy class rel. boundary of $\gamma$ and additive w.r. to juxtaposition of paths.  Indeed if $\gamma[a,b]\subset U_\alpha $ for some $\alpha$ then $\int_{\gamma} \omega= f_\alpha(b)- f_\alpha(a).$ If not, one choses a subdivision of $[a,b],$  $a= t_0 <t_1 <\cdots t_r= b,$ such that $\gamma_i:= \gamma|_{ [t_i, t_{i+1}]}$ lie in some open set $U_\alpha$ and put $\int_{\gamma} \omega:= \sum \int_{\gamma_i} \omega.$ This assignment  defines an homomorphism $\xi(\omega): H_1(X;\mathbb Z)\to \mathbb R,$ equivalently a cohomology class $\xi(\omega).$}
\vskip .1in 
One denotes by 
{$\mathcal Z^1(X)$ the set of all TC1-forms }and by
{$\mathcal Z^1(X; \xi):= \{\omega\in \mathcal Z^1(X)\mid  \xi(\omega)= \xi\}.$}
\vskip .1in

 Let 
$\xi\in H^1(X;\mathbb R)= Hom (H_1(X;\mathbb Z), \mathbb R)$ and  $\Gamma = \Gamma (\xi):
= (img \xi) \subset \mathbb R.$  If $X$ is a compact ANR then  $\Gamma\simeq \mathbb Z^k$ and one refers to $k$ as  the {\it degree of irrationality} of $\xi.$
\vskip .1in

The surjective homomorphism {$\xi: H_1(X;\mathbb Z)\to \Gamma$} defines the associated $\Gamma-$principal cover, {$\pi: \tilde X\to X$}  i.e. the free action {$\mu:\Gamma\times \tilde X\to \tilde X$} with $\pi $ the quotient map $\tilde X \to \tilde X/\Gamma= X.$
\vskip .1in

 Consider 
$f: \tilde X\to \mathbb R,$ $\Gamma-$equivariant maps, i.e. satisfying $\boxed{f(\mu(g,x))= f(x)+g}.$ 

\vskip .1in
{\bf Second definition}
\begin{definition}  \

A   TC1- form $\omega$ of cohomology class $\xi$ is an equivalence class of $\Gamma-$equivariant maps  $f:\tilde X\to \mathbb R$ with $f_1\sim f_2$ iff $f_1-f_2$ is locally constant.
\end{definition}
Such map $f$ is referred to as a {\it lift} or a representative of  $\omega$  and one writes 
$f\in \omega$, $\omega\in [\omega]=\xi(\omega)$ where  $\xi(\omega)= \xi.$  

 One denotes by 
{$\mathcal Z^1(X;\xi)$} the set of TC-1 forms in the class $\xi$  and by 
{$\mathcal Z^1(X)= \bigcup_{\xi\in H^1(X;\mathbb R)} \mathcal Z(X; \xi).$}
\vskip .2in

Clearly the two definitions, first and second, are equivalent.

\vskip .2in
 {\bf Tameness for TC1-forms}
\begin{definition}\

A continuous map  $f:\tilde X\to \mathbb R$ is called {\it weakly tame} if the following hold: 

(i) For any $I\subseteq \mathbb R$ closed interval $f^{-1}(I)$ is an ANR.   

(ii) For $t\in \mathbb R,$ $R^f(t):= \dim H_r(\tilde X^f_{t}, \tilde X^f_{< t} ) + \dim H_r(\tilde X_f^{t}, \tilde X_f^{>t})<\infty$  

where $\tilde X^f_{t}= f^{-1}((-\infty, t]), \ \tilde X^f_{< t}= f^{-1}((-\infty, t))$  and $\tilde X_f^{t} = f^{-1}([t,\infty)), \ \tilde X_f^{>t}= f^{-1}((t,\infty)).$ 

The value $t\in \mathbb R$ is called {\it regular} if $R^f_t=0$ and {\it critical} otherwise.

(iii) The set of critical values $CR(f)$ is at most countable.  
\end{definition}
 When there is no specification of coefficients in the notation  it is understood that the homology $H_r(\cdots)$  is with coefficients in the field $\kappa$ and $H_r(\cdots)$ is a $\kappa-$vector space. 

Note that if $f:\tilde X\to \mathbb R$ is $\Gamma-$equivariant then $CR(f)$ is $\Gamma-$invariant.
\begin{definition}\

Let $X$ be a compact ANR.
The TC1- form $\omega\in \mathcal Z^1(X:\xi),$ 
is called {\it tame} if one lift 
$f:\tilde X\to \mathbb R$ (and then any other lift) is weakly tame and the set
$CR(f)/\Gamma$ is finite.
 \end{definition}  

One denotes by 
$\mathcal Z^1_{tame}(X;\xi)$  the subset of $\mathcal Z^1(X;\xi)$ consisting of tame   TC-1 forms
by and $\mathcal Z^1_{tame}(X):=  \bigcup_{\xi\in H^1(X;\mathbb R)} \mathcal Z_{tame}(X; \xi).$  

\vskip .2in 

 {\bf Topology on the set $\mathcal Z^1(X:\xi)$ } 

Suppose $X$ is a compact ANR. The set $\mathcal Z^1(X:\xi)$ can be equipped with the complete metric $D(\omega_1, \omega_2)$ (whose induced topology is referred to as the {\it compact open topology}) defined as follows:
$$D(\omega_1, \omega_2)= \inf _{\tiny \begin{aligned} f_1\in \omega_1\\f_2\in \omega_2\end{aligned}} D(f_1, f_2)$$ with $D(f_1, f_2):= \sup _{x\in \tilde X} |f_1(x)- f_2(x)|.$

\section {The set 
$\mathcal O(\omega)$ and spaces of configurations} 

As  pointed out, if $f:\tilde X\to \mathbb R$ is $\Gamma-$invariant then the set $CR(f)$ of critical values is $\Gamma-$invariant. Denote  
$\mathcal O(f):=CR(f)/\Gamma.$ 
If $f_1$ and $f_2$ are two such lifts  and $X$ is connected  (otherwise the considerations are done component-wise), hence $f_2= f_1 +s,$ then the translation by $s$ gives an identification $\theta(s) :\mathcal O(f_1)\to \mathcal O(f_2).$

The set $\mathcal O(\omega)$, of {\it orbits of critical values} is defined by
 $\mathcal O(\omega):= \sqcup_{f \in \omega} \mathcal O(f)/ \sim$ with $o\sim o',$ $o\in \mathcal O(f_1)$ $o'\in \mathcal O(f_2),$ iff $f_2= f_1 +s$ and 
 $\theta(s)(o_1)= o_2.$
\vskip .2in
{\bf Configurations}

A map $\delta: Y\to \mathbb Z_{\geq 0},$ $Y$ topological space,  $\delta$ a map with finite support is called {\it configuration}. One defines  $\supp \ \delta:=\{x\in Y\mid \delta(x)\ne 0\}$ and $\sharp \supp \ \delta:= \sum_{x\in Y} \delta(x).$ 

Denote by $Conf (Y)$ the set of all configurations on $Y$ and by $Conf_N (Y)$ the subset  $Conf_N (Y):= \{\delta \in Conf(Y) \mid \sharp \supp\ \delta= N\}.$  Since $Y$ is a topological space $Conf_N(Y)= Y^N/\Sigma_N,$  
$\Sigma_N$ the permutaion group of $N$ elements, is equipped with the obvious topology, the same as the {\it collision topology}  and when $Y= Z\setminus K, K $ closed subset of $Z$, $Conf(Y)$ is  equipped with the {\it bottleneck topology} w.r. to the pair $(Z, K)$. 

A fundamental system of neighborhoods for a configuration $\delta \in Conf _N(Y)$ with support $\{y_1, y_2, \cdots y_r\}$ in the collision topology is provided by the sets of configurations $\mathcal  U(\delta, \{U_i\}):= \{\delta' \in Conf (Y) \mid \sum_{y\in U_i}  \delta'(y)= \delta(y_i)\}$ with $U_i$ a disjoint collection of open neighborhoods of $y_i.$ 

A fundamental system of neighborhoods for a configuration $\delta \in Conf (Y)$ with support $\{y_1, y_2, \cdots y_r\}$ in the bottleneck topology is provided by the sets $\mathcal  U(\delta, \{U_i\}, V):= \{\delta'\in Conf (Y) \mid \sum_{y\in U_i}  \delta'(y)= \delta(y_i), \supp\ \delta' \subset (V\sqcup \bigcup_{i=1, \cdots r}U_i )\}$ with $U_i$ and $V$ a disjoint collection of open sets of $Y,$  $U_i$ neighborhood of $y_i$  and $V$ neighborhood of $K.$ 

\section {The Alternative  Novikov theory  for topological closed one  form of arbitrary degree of irrationality}\

Denote by $\mathcal Z^1_{tame} (\cdots):= \{ \omega \in \mathcal Z^1 (\cdots) \mid \omega \
 \text{tame}\}.$
In what follows 
for any tame  TC-1 form $\omega \in \mathcal Z(X;\xi),$ field $\kappa$ and nonnegative integer $r$ one defines two configurations 
$$\boxed{\underline \delta^\omega_r: \mathbb R\to \mathbb Z_{\geq 0}}\ \text{and} \  \boxed {\underline \gamma^\omega_r: \mathbb R_+\to \mathbb Z_{\geq 0}}.$$ 
The blue circles in Figure 1. represent the closed $r-$bar codes located at $\mathbb R \ni x \leq 0$ and the red circles represent the open $(r-1)-$barcodes located at $\mathbb R\ni x <0.$
The brown circles in Figure 2. represent the closed-open $r-$barcodes located on  $\mathbb R_+= (0,\infty).$ 

\begin{figure}\ 

\begin{center}
\begin{tikzpicture}
\draw [ultra thick]  (-6,0) -- (6,0);
\draw [dashed, ultra thick]  ((6,0) -- (7,0);
\draw [dashed, ultra thick]  ((-7,0) -- (-6,0);
\draw[fill=blue] (0,0) circle [radius=0.1cm];
\draw[fill=blue] (2,0) circle [radius=0.1cm];
\draw[fill=blue] (4,0) circle [radius=0.1cm];

\draw[fill=red] (-1,0) circle [radius=0.1cm];
\draw[fill=red] (-2.5,0) circle [radius=0.1cm];
\draw[fill=red] (-4,0) circle [radius=0.1cm];
\end{tikzpicture}
\end{center}
\caption {Configuration 
{{$\underline\rho^\omega_r$}}}
\label{FIG52}
\end{figure}
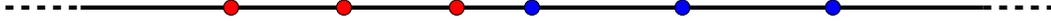

\vskip .4in 
\begin{figure}\ 

\begin{center}
\begin{tikzpicture}
\draw[ultra thick, <-]  (9.6,0) -- (16,0);
\draw[dashed, ultra thick]  ((16,0) -- (17,0);
\draw[fill=brown] (11.6,0) circle [radius=0.1cm];
\draw[fill=brown] (13.6,0) circle [radius=0.1cm];
\end{tikzpicture}
\end{center}
\caption {Configuration {{$\underline \gamma^\omega_r$}}}
\label{FIG52}
\end{figure}
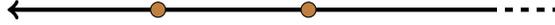

Define also the numbers 
 $\beta^\omega_r, \rho^\omega_r \ \text{and} \ c^\omega_r,$
\vskip .2in  
\hskip 1in $\boxed{\beta_r^\omega:=\sum_{t\in \mathbb R} \underline \delta_r^\omega (t)},$\quad  
$\boxed{\rho_r^\omega:=\sum_{t\in \mathbb R_+} \underline \gamma_r^\omega (t)}$ 
\vskip .2in
and \hskip 1in  $\boxed{c_r^\omega= \beta^\omega_r + \rho^\omega_r + \rho^\omega_{r-1}} .$
\vskip .2in

These configurations satisfy the Theorems \ref{TT}, \ref{TP} and \ref{TS} below.

\vskip .1in
{\bf The results}
\vskip .1in

\begin{theorem} (Topological properties)        \label {TT}\

Let $\omega \in \mathcal Z^1_{tame}(X) .$  One has: 
\begin {enumerate} 
\item $ \underline \delta^\omega_r(t)\ne 0$ or $\underline \gamma^\omega_r (t)\ne 0$  
implies that for any $f\in \omega,$  $t= c' - c''$ with $c', c''\in CR(f),$
\item $$\beta_r^\omega=\boxed{\sum _{t\in \mathbb R}\underline \delta^\omega_r(t)= \beta^N_r(X;\xi)}$$
\item for any $f\in \omega$ and $a_o\in o\in \mathcal O(f),$ 
$$ c^\omega_r=  \boxed {\begin{aligned} \sum_{t\in \mathbb R}& \underline \delta^\omega_r(t)+
\sum_{t\in \mathbb R_+} \underline \gamma^\omega_r(t)+ \sum_{t\in \mathbb R_+} \underline \gamma^\omega_{r-1}(t) = \\
= &\sum _{o\in \mathcal O(f)} \dim H_r(\tilde X^f_{a_o}, \tilde X^f_{<a_o})\end{aligned}}$$
{where  $ \tilde X^f_{a_o}= f^{-1}((-\infty, a_0]) , \tilde X^f_{<a_o}= f^{-1}((-\infty, a_0))$.}
\end{enumerate}
\end{theorem}

Note that  
if $X$ is a smooth closed manifold and $\omega$ is a smooth Morse form Lyapunov for a Morse-Smale vector field then the number of zeros of $\omega$  (= the rest points of the vector field) of Morse index $r$  
is equal to 
$ \sum_{o\in \mathcal O(f)} \dim H_r(\tilde X^f_{a_o}, \tilde X^f_{< a_o}).$ 

\begin{corollary}
The Novikov complex of the Morse form $\omega$ in Hodge form  is  given by  $C_r= C_r^-  \oplus  C_r^+  \oplus  \mathcal H_r$
with 

$\mathcal H_r= \kappa^{\beta_r^\omega},$

$C_r^-= C_{r-1}^+=\kappa^{\rho_r^\omega},$

$\partial_r= \begin {pmatrix} 0&0&0\\Id&0&0\\0&0&0
\end{pmatrix}.$ As in \cite {B} in the case of an angle-valued map, this complex $(C\ast, \partial_\ast)$ in Hodge form can be referred to as  the {\it AN complex associated to $\omega.$}
\end{corollary}
\vskip .2in
 {\it Relation to other work \cite{UZ}and \cite{BH}}

The points in $\supp \underline \delta^\omega_r\sqcup \supp \gamma^\omega_r \sqcup \supp \gamma^\omega_{r-1}$ with their multiplicity correspond bijectively to  UZ (=Usher-Zhang) verbose  $r-$bar codes.  
and in the case $\xi(\omega)$ integral to the BD(=Burghelea-Dey)- closed $r-$barcodes plus BD-open $(r-1)-$barcodes  plus closed-open $r-$barcodes plus closed-open $(r-1)-$barcodes cf \cite{BD11} or \cite {BH}.
In the case $\xi(\omega)$ is trivial the points in $\supp \underline \gamma^\omega_r$ with their multiplicity correspond bijectively to the finite barcodes 
in classical ELZ-persistence (cf. \cite {ELZ}) and  in the case $\xi(\omega)$ is integral to the closed-open $r-$barcodes in angle valued persistence.
\vskip .2in

\begin{theorem} (Poincar\'e duality properties) \label {TP}\

Suppose $X$ is a closed topological $n-$dimensional manifold and $\omega\in \mathcal Z^1_{tame}(M).$ Then
\begin {enumerate} 
\item $ \underline \delta^\omega_r(t)= \underline \delta^\omega_{n-r}(-t),$ 
\item $ \underline \gamma^\omega_r(t)= \underline \gamma^{-\omega}_{n-r-1}(t).$  
\end{enumerate}
\end{theorem}

Recall that:

-  The set of TC-1 forms $\mathcal Z^1(X;\xi)$ is an $\mathbb R-$vector space equipped with the norm given by the distance $D(\omega_1, \omega_2)$ defined above. 
 
-  The space of configurations $Conf_N(Y):= \{ \delta: Y\to \mathbb Z_{\geq 0}\mid \sum \delta(y)= N\}$ is equipped with the obvious {\it collision topology}.
 
-  For $K$ closed subset of $Y$  the space of configurations 
 $Conf (Y\setminus K)$ is equipped with the {\it bottleneck topology} provided by the pair $ K \subset Y.$
 With respect to these topology one has 
\begin{theorem} (stability properties) \label {TS}\

\begin {enumerate} 
\item  The assignment $\mathcal Z^1_{tame}(X;\xi)\ni \omega \rightsquigarrow \underline \delta^\omega_r \in Conf_{\beta^N_r(X;\xi)}(\mathbb R)$  is continuous.
\item  The assignment $\mathcal Z^1_{tame}(X;\xi)\ni \omega \rightsquigarrow \underline \gamma^\omega_r\in Conf( [0,\infty)\setminus 0)$  is continuous.
\end{enumerate}
\end{theorem} 
\vskip .1in

{\bf Description of the configurations $\underline \delta_r^\omega$ and $\underline \gamma_r^\omega.$}
\vskip .1in
Fix a field $\kappa.$
Let $\pi :\tilde X\to X$ be a $\Gamma-$principal covering  with the free action $\mu :\Gamma\times \tilde X\to \tilde X$  where $\Gamma\subset \mathbb R$  is a f.g. subgroup of $\mathbb R.$  Let $f: \tilde X\to \mathbb R$ be a $\Gamma-$equivariant map, i.e. $f(\mu(g,x))= f(x)+g$ which is weakly tame and $\mathcal O(f)$ finite.

Recall that  $\tilde X^f_a = f^{-1}((-\infty ,a]), \tilde X^f_{<a}= f^{-1}((-\infty ,a)), \tilde X_f^a= f^{-1}([a,\infty)),  \tilde X_f^{>a}= f^{-1}((a, \infty))$
and denote 
\begin{enumerate}
\item 
{$\boxed{\mathbb I^f_a(r): = img (H_r(\tilde X_{a}) \to H_r(\tilde X))}, $} \quad    $\mathbb I^f_a(r)\subset H_r(\tilde X)$
\vskip .1in
\hskip .2in $\mathbb I^f_{<a}(r):= \cup_{a'<a} \mathbb I^f_{a'}(r),$ \quad 
$\mathbb I^f_{-\infty}(r)= \cap_a \mathbb I^f_{a}(r),$
\vskip .1in

\item $\boxed{\mathbb I_f^a(r): = img (H_r(\tilde X^{ a})  \to H_r(\tilde X))}, $  \quad  $\mathbb I_f^a(r)\subset H_r(\tilde X)$
\vskip .1in
\hskip .2in $\mathbb I_f^{>a}(r):= \cup_{a'>a} \mathbb I_f^{a'}(r)$, \quad $\mathbb I_f^{+\infty}(r)= \cap_a \mathbb I_f^{a}(r),$ 
\vskip .1in

\item 
For $a<b$ \ \ 
{$\boxed{\mathbb T^f_r(a,b): = \ker (H_r(\tilde X_{a})  \to H_r(X_{b}))}, $   
\vskip .1in
\hskip .2in $\mathbb T_r(a, <b):= \cup_{b'<b} \mathbb T_r(a, b')  \subseteq \mathbb T_r(a,b),$ 
\vskip .1in 
\item For $a'<a$\   let  $i: \mathbb T_r(a',b)\to \mathbb T_r(a,b)$ the induced map 
and  
\vskip .1in
\hskip .2in $\mathbb T_r(<a,b):= \cup_{a'<a} i (\mathbb T_r(a',b)) \subseteq \mathbb T_r(a,b).$}
\end{enumerate} 
 
Define  
\vskip .1in
for $a,b\in \mathbb R$ 
{$$\boxed{\hat\delta^f_r(a,b):= \frac {\mathbb I_a(r)\cap \mathbb I^b(r)}{\mathbb I_{<a}(r)\cap \mathbb I^b(r) + \mathbb I_a(r)\cap \mathbb I^{>b}(r)},\   \delta^f_r(a,b):=\dim \hat\delta^f_r(a,b),} $$}

 for $a,b\in \mathbb R,  a<b$
{$$\boxed{\hat\gamma^f_r(a,b):= \frac {\mathbb T_r(a,b)}{\mathbb T_r(<a, b)+  T_r(a, <b)} ,\ \    \gamma^f_r(a,b):=\dim \hat\gamma^f_r(a,b).}$$}

and for $t\in \mathbb R$ 
{$$\boxed{\mathbb F^f_r(t):= \sum_ {\tiny \begin{cases}a-b \leq t, \\
a,b\in CR(f)\end{cases}} \mathbb I_a(r)\cap \mathbb I^b(r)}$$}

One can show  that: 

\begin{proposition}\

\begin{enumerate}
\item {$\delta^f_r(a,b)= \delta^f_r(a+g, b+g)$} and 
{$\gamma^f_r(a,b)= \gamma^f_r(a+g, b+g)$} for any $g\in \Gamma,$

\item for $a$ regular value {$\supp \delta^f_r \cap (a \times \mathbb R)$} and {$\supp \gamma^f_r\cap (a\times (a,\infty))$} is empty,
\item for $a$ critical value {$\supp \delta^f_r\cap (a\times \mathbb R)$} and  {$\supp \gamma^f_r\cap (a\times (a,\infty))$} is finite,

\item for $b$ regular value {$\supp \delta^f_r\cap (\mathbb R \times b)$} and {$\supp \gamma^f_r\cap ((-\infty, b)\times b)$} is empty,
\item for $b$ critical value {$\supp \delta^f_r\cap (\mathbb R\times b)$} and  {$\supp \gamma^f_r\cap ((-\infty, b)\times b)$} is finite.
\end{enumerate}
\end{proposition}

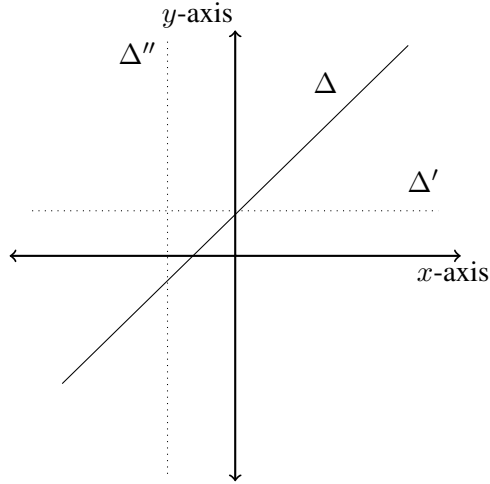
\begin{figure}
\begin{center}
\begin{tikzpicture}
\draw [thick,<->]  (0,3) -- (0,0) -- (3,0);
\node at (-0.5,3.2) {$y$-axis};
\node at (2.9,-0.2) {$x$-axis};
\node at (-1.3,2.7) {$\Delta''$};
\node at (1.2,2.3) {$\Delta$};
\node at (2.5,1) {$\Delta'$};
\draw [thick,<->]  (0,-3) -- (0,0) -- (-3,0);
\draw (-2.3,-1.7) -- (2.3,2.8);
\draw [dotted] (-2.7,0.6) -- (2.7,0.6);
\draw [dotted] (-0.9,-2.9) -- (-0.9,2.9);
\end{tikzpicture}
\end{center}
\caption {the domain of $\delta^f_r$}
\label{FIG52}
\end{figure}

Note that :

\begin{enumerate}
\item
For any 
line $\Delta'$ or $\Delta''$ parallel to either 
$x-$axis or $y-$axis the cardinality of $\supp \delta^f_r\cap \Delta$ or $\supp \gamma^f_r\cap \Delta$ is finite.
\item  For any line $\Delta$ parallel to the  first diagonal $\{(x,y)\in \mathbb R^2 \mid x-y=0\}$ the values of each $\delta^f_r$ and $\gamma^f_r$ are the same on each point of $\supp \delta^f_r\cap \Delta$ resp. of $\supp \gamma_r^f\cap \Delta.$
\item Suppose the line $\Delta$ is parallel to the  first diagonal $\{(x,y)\in \mathbb R^2 \mid x-y=0\}$ and  $\omega$ is of  degree of irrationality $>1.$ If  $\supp \delta^f_r\cap \Delta$ and  $\supp \gamma^f_r\cap \Delta$ is nonempty then it is  everywhere dense in $\Delta.$
\end{enumerate}
\vskip .1in 

\begin{obs}
For $a_o\in o\in \mathcal O(f),$ $f\in \omega,$ 
{$\delta^f_r(a_o,a_o +t)$}  and 
$\gamma^f_r(a_o,a_o +t)$\ 
are independent on the choice of $a_o$ in $o$ and of  $f$ in $\omega.$ 
\end{obs}

Define 
{$\boxed{\underline \delta^\omega_{o,r} (t):= \delta^f_r(a_o,a_o +t)}$}
and then 
{$$\boxed{\underline \delta_r^\omega(t)= \sum _{o\in \mathcal O(\omega)} \underline \delta^\omega_{o,r} (t)}.$$}
\vskip .2in

Define 
{$\boxed{\underline \gamma^\omega_{o,r} (t):= \gamma^f_r(a_o,a_o +t)}$}
and then 
$$\boxed{\underline \gamma_r^f(t)= \sum _{o\in \mathcal O(\omega)} \underline \gamma^f_{o,r} (t)}.$$

\section{About the proof of Theorems \ref{TT},\ \ref{TP}  and\  \ref{TS}}

\begin{proposition}\label{P51}\
Let $f: \tilde X\to \mathbb R,$ with $f\in \omega\in \mathbb Z^1_{tame}(X).$  
\begin{enumerate}
\item For $t \in \mathbb R$  one has 
$$H_r(\tilde X^f_{t}, \tilde X^f_{<t})\simeq  {\mathbb I^f_a(r)/ \mathbb I^f_{<a}(r)} \oplus {\coker  (\mathbb T^f_r(<t,\infty)\to  \mathbb T^f_r(t, \infty))} \oplus  \ker  (H_{r-1}(\tilde X^f_{< t})\to  H_{r-1}(\tilde X^f_{t})).$$
and if $t$ is a regular value then $H_r(\tilde X^f_{t}, \tilde X^f_{<t})=0.$
 
\item ${\mathbb I^f_{ t}(r) / \mathbb I^f_{<t}(r)\simeq \oplus _{s\in \mathbb R} \hat \delta^{f}_r(t,s)}.$

\item ${\coker  (\mathbb T^f_r(< t,\infty)\to  \mathbb T^f_r(t, \infty)) \simeq \oplus _{s\in \mathbb R_+} \hat \gamma^{f}_r(t, t+s)}.$

\item ${\ker  (H_{r-1}(\tilde X^f_{< t})\to  H_{r-1}(\tilde X^f_{\leq t})) \simeq \oplus _{s\in \mathbb R_+} \hat \gamma^{f}_{r-1}(t-s, t)}.$
\end{enumerate}

\end{proposition} 

The proof of item 1  follows from the inspection of the homology long exact sequence of the pair $(\tilde X_{\leq t}, \tilde X_{<a}).$   The proof of the other items 
follow from the finite dimensionality of $H_r(\tilde X_{\leq t}, \tilde X_{<a})$ and the definitions of $\delta^f_r(a,b)$ and $\gamma^f_r(a,b).$

The results  remain true for any $f: \tilde X\to \mathbb R$ weakly tame and not necessary  $\Gamma-$equivariant, in particular for any restriction of $f$ to an open subset $U\subset \tilde X.$    
\vskip .2in

In view of the action of $\Gamma$ on $\tilde X$  the $\kappa-$vector spaces  
  $H_r(\tilde X),$ \  $\mathbb I^f_{-\infty}(r), \mathbb I_f^\infty(r) \subseteq H_r(\tilde X)$ are actually f.g. $\kappa[\Gamma]-$module resp. submodules. 
Denote by   $TH_r(\tilde X)$ the submodule of torsion elements of $H_r(\tilde X).$
\vskip .1in

 For $f\in \omega\in \mathcal Z^1_{tame}(X)$  consider the filtration of f.g. $\kappa[\Gamma]-$modules indexed by $t\in \mathbb R$

$$\mathbb I_{-\infty}(r)= \mathbb F^f_r(-\infty)\subseteq \cdots  \mathbb F^f_r(t')\subseteq \mathbb F^f_r(t) \cdots \subseteq \mathbb F^f_r(\infty)= H_r(\tilde X),  t'<t $$ 
and observe that in view of f,g. of $H_r(\tilde X)$  there are finitely many $t'$s  where  the rank of $\mathbb F_r(t)$ changes.
\vskip .2in 

Theorem \ref {TT} (Topological properties) follows  from  Proposition \ref{P52} below and  the following facts: 

-    $\mathbb F^f_r(t)/ TH_r(\tilde X)\simeq \oplus_ {\tiny \begin{cases}a-b \leq t, \\
a,b\in CR(f)\end{cases}}  \hat \delta^f_r(a,b) ,$

-   if $a$ is a regular value then $H_r(\tilde X^f_{a}, \tilde X^f_{<a})= 0$  and if $a\in o\in \mathcal O(f)$ then 
$$\dim H_r(\tilde X_{ a}, \tilde X_{<a})= \sum_{t\in \mathbb R}\underline \delta^\omega_{o,r} (t)+ \sum_{t\in \mathbb R_+}\underline \gamma^\omega_{o,r}(t) +\sum _{t\in \mathbb R_+}\underline \gamma^\omega_{o,r-1}(t),$$
which follow from  Proposition \ref{P51} above. 
\begin{proposition}\label {P52}\

$\mathbb I^{f}_{-\infty} (r)= \mathbb I^{\infty}_{f}(r)= TH_r(\tilde X)$
\end{proposition}
The proof of this proposition  is similar to the proof of Proposition 5.10 in \cite{B}.
\vskip .2in
 
To prove Theorem \ref{TP} 
one uses Borel-Moore homology. One  defines the analogues $^{BM} \hat \delta_r$ and $^{BM} \hat \gamma_r$ of  $\hat \delta_r$ and $\hat \gamma_r$ and then  $^{BM} {\hat \delta}_r$ and  $^{BM} \hat \gamma_r.$
and one proves  
\begin{proposition}\label {P53}\
$^{BM}\delta ^{f}_r(a,b)= \delta^{f}_r(a,b)$ and \ \ $^{BM}\gamma ^{f}_r(a,b)= \gamma^{f}_r(a,b).$
\end{proposition} 
  One derives these equalities  using similar arguments as in the proof of Theorems 5.6 and 6.2 in \cite {B} (which contain Theorem \ref{TP} as a particular case  when $\omega$ of degree of irrationality zero or one). However the arguments are longer since $f$ is not proper in general.

\vskip .2in
To prove item 1 in Theorem \ref {TS} one verifies first 
\begin{proposition}\label {P54} 
Suppose $f\in \omega_1, g\in \omega_2,$ \ \ $\omega_1, \omega_2 \in \mathcal Z_{tame}(X;\xi).$ If  $|f-g|<\epsilon$ then 
$$\mathbb F^f_r(t) \subseteq \mathbb F^g_r(t+2\epsilon) \subseteq \mathbb F^f_r(t+4\epsilon)$$  
\end{proposition}

Once Proposition \ref{P54} is established one proceeds as in  \cite {B},  Proposition 5.7 and Theorem 5.2. 

To check  item 2. one proceeds as follows.

\begin{itemize}
\item Consider the configuration $\lambda^\omega_r: \mathbb R_+\to \mathbb Z_{\geq 0}$ defined  by $$\lambda^\omega_r(t)= \underline \delta^\omega_r(t)+ \underline \gamma_r^\omega(r)+ \underline \gamma^{\omega}_{r-1}(t),$$ and establish the continuity of the assignment 

$\mathcal Z_{tame}(X;\xi)\ni \omega \rightsquigarrow \lambda_r^{\omega}\in Conf([0,\infty)\setminus 0)$ as in the proof of the CEH-stability cf.  \cite{CEH07}. 

\item Use item 1 to conclude  the continuity for the assignment given by  $\delta^\omega_{r, +}$  the restriction of the configuration  $\underline \delta ^\omega_r $  to $\mathbb R_+.$
\end{itemize}

Clearly the continuity of the assignment $\omega\rightsquigarrow \underline \gamma_r^\omega(r)+ \underline \gamma^{\omega}_{r-1}(t)$ and then of $\omega\rightsquigarrow\underline \gamma_r^\omega(r)$ is the same as   the continuity of  the assignment 
$\omega\rightsquigarrow\lambda^\omega_r(t) -\underline \delta^\omega_{r, +}(t).$
The continuity of  $\lambda^\omega_r(t) $ follows also from  \cite {UZ}.

\section {Comments and Applications}

The  results summarized above  can  and will be used for additional explorations / results in: 

1. Dynamics,

2. Topology of smooth manifolds with symmetry,

3. Topology of the free loop space, 

3. Geometrization of Data.
 
\vskip .1in
\noindent {\bf Dynamics}: One can extend the results of Novikov theory from smooth flows on compact manifold with a Morse closed  one form as Lyapunov to 
dynamics defined by a continuous flow on a compact ANR whose trajectories minimize an "action".  Precisely one can detect  rest points  instantons and closed trajectories from the topology  of the underlying space 
as in the case of classical Morse-Novikov theory. A paper on such applications is in preparation.

\vskip .1in
\noindent{ \bf Morse-Novikov theory on a G-Manifolds}: 
Given a  smooth manifold $M,$  a closed one form and a smooth vector field with the closed form Lyapunov for the vector field which are invariant  to a smooth action of a compact Lie group, by passing to the "quotient space", one obtains a compact ANR, a tame TC1-form and a continuous flow on the ANR. It seems interesting to compare the barcodes of the TC1-form with the complex which calculates the $G-$Novikov homology derived via  classical  G-Novikov theory. 
 This will be done in a future work. 

\vskip .1in
\noindent{\bf Morse-Novikov theory for free loop space}:  Suppose $M$ a closed Riemannian  manifold  and $M^{S^1}$ the infinite dimensional manifold of smooth free loops on $M$. For a fixed integer $N$ and a fixed $\epsilon >0$ denote by $M^N_\epsilon:= \{(x_1, x_2, \cdots x_n)\mid d(x_i, x_j)\leq \epsilon, d(x_n, x_1)\leq \epsilon\}.$ Clearly $M^N_\epsilon$ is a smooth compact manifold with corners hence a compact ANR.  If $\epsilon$ is small enough  to insure unicity of a geodesic between two  points $x,y$ with $d(x,y) \leq \epsilon$  then one has an obvious inclusion $M^N_\epsilon \subset M^{2N}_{\epsilon/2} \subset M^{S^1}.$  A smooth closed one form $\omega$ on $M^{S^1}$ induces the smooth closed one form $\omega_{N,\epsilon}$ on $M^n_\epsilon.$  Clearly  for any  $r$ there exists an integer  $K(r)$ s.t. for $k> K(r)$ the inclusion $M^{2^k n}_{\epsilon/ 2^k}\subset M^{S^1}$ induces an isomorphism in $r-$ dimensional homology and cohomology.   One expects that the barcodes  of the pair $M^{2^k n}_{\epsilon/ 2^k},\omega|_{M^{2^k n}_{\epsilon/ 2^k}}$ stabilize with $k$ large enough and permit the definition of  barcodes for $(M^{S^1},\omega).$  We plan to apply this to  a system $(M,\sigma, g, J)$  with $(M.\sigma)$ a symplectic manifold equipped with a compatible almost complex structure $J$ (the Riemannian metric being determined by $(\omega, J)$) and a  time dependent periodic symplectic vector field, cf \cite {BH2}. These data provides a closed one form on $M^{S^1}$ determined by $\sigma$ and  the closed one form induced from  a  time dependent periodic symplectic vector field, cf \cite {BH2}. This is of interest in symplectic topology and Hamiltonian dynamics. This is work in progress.       

\vskip .1in
\noindent {\bf Geometrization of some data}:  Suppose $(X,d)$ is a finite metric space (i.e. point cloud data), $S\subset X\times X$ and $f: S\to \mathbb R$ with the property that $(x,y)\in S$ implies  $(y,x)\in S$ and $f(x,y)= -f(y,x).$
Let $\epsilon (f,d)>0$ be  the largest real number $\epsilon $ such that 
for any three points $x,y,z\in X$ with $(x,y), (y,z), (z,x) \in S$ and $d(x,y), d(x,z), d(y,z) <\epsilon$ one has  $f(x,y)+ f(y,z) + f(z,x)=0.$   

Clearly the Rips complex for any $\epsilon <\epsilon (f,d)$ is a simplicial complex equipped with a TC1-form $\omega(\epsilon)$ given  by the simplicial cocycle  defined by $f.$  
It might be interesting to study the bar codes based on  persistent (Novikov) homology or even  of the family of the $AN-$chain complexes  associated with $\omega(\epsilon),$ a considerable refinements of Novikov homology. They offer additional topological invariants to describe features of $\{(X, d), f\}.$

\end{document}